\documentclass{article}%
\usepackage{amsmath}%
\usepackage{amsfonts}%
\usepackage{amssymb}%
\usepackage{graphicx}%
\usepackage{color} 

\def\l{\lambda}

\def\s{\sigma}

\def\Bbb{\mathbb}

\def\diag{\mathop{\rm diag}}

\def\T{\mathop{\footnotesize{t}}}

\begin{document}
\title{}
\author{}
\date{}
\maketitle

\noindent
Book Review{\footnote{to appear in Linear Algebra Appl.}}
\bigskip

\noindent
{\bf Matrix Theory,  by Xingzhi Zhan,
American Mathematical Society (AMS), Graduate Studies in Mathematics (GSM) Volume 147,
2013,  x+253~pp., Hardcover,
 ISBN 978-0-8218-9491-0}

\bigskip
Zhan's book {\em Matrix Theory} is a masterpiece in the field.
It presents  perspectives of modern matrix theory: fundamentals, new developments, open problems, applications, as well as connections
with other mathematical subjects. It shows that matrix theory is a rich research branch of mathematics.

This is a book written by a matrix expert and it is  pleasant to read.

Matrices play an important role in all branches of mathematics, applied and pure; the theory of matrices is the core of linear algebra,
while linear algebra is the foundation of all mathematics. Of the about 150 books in the GSM series published by AMS, this is the only one devoted to matrix theory. It deserves attention.

Three brief reviews of the book have appeared; see \cite{Berg, DFarenickRev, RJun}. I would like to give a more detailed account, in hope that
 it may help   instructors   teach a graduate
 course on matrix theory when using the book as a text.

It is my personal opinion that the general mathematical community misinterprets {\em matrix theory} as {\em linear algebra}, and vice versa. Matrix theory is not identical to linear algebra, because most modern matrix problems and the research methods and techniques for the matrix problems are neither linear nor algebraic. It appears to this reviewer that analysis and combinatorics have played more important roles than algebra and other subjects in  modern matrix theory research. But matrix theory and  linear algebra are two sister subjects of mathematics after all.
I think this book confirms my point of view.

 This book selects many important results with  simple and elegant proofs. Varying in depth and variety, these results  are central to matrix theory.
 Selected topics are  attractive
 and exciting and show the beauty of the field. There are a plenty of exercises in most of the chapters.  The book
  may be used as a standard two semester  course at the graduate level.

Prerequisite is a solid background in linear algebra 
and analysis. Occasionally, an ``advanced" terminology (such as {\em Gr\"{o}bner basis}, p.~26;
{\em homotopy}, p.~105; {\em Borel measure}, p.~146) is used; 
but it can be easily found in a text and should be accessible to graduate students without difficulty.

Let me begin with a few easy to understand but hard to answer questions. They are
 easy because every undergraduate student can understand them and he/she may even ask the questions during a class; these questions are hard to answer since
  advanced tools are needed in their solutions. The answers to these questions can be found in the book.

{\bf ``Easy" questions:}

\begin{enumerate}
\item Is every complex number an eigenvalue of some nonnegative matrix? 
That is: Let $\l$ be a complex number. Is it possible that $Ax=\l x$ for some nonnegative matrix $A$ and a column vector $x\neq 0$?
(See p.~131.)

\item How many negative entries can the square of a real matrix have?
(See p.~131.)

\item If two real matrices are (unitarily) similar over the complex number field $\Bbb C$, are they (orthogonally) similar over the real number field $\Bbb R$?
That is: Let $A$ and $B$ be $n\times n$ real matrices. If $U$ is an $n\times n$  nonsingular complex (respectively,  unitary) matrix such that $U^{-1}AU=B$, is it possible
that $V^{-1}AV=B$ for some $n\times n$ nonsingular real (respectively,  orthogonal) matrix $V$?
(See pp.~181--182.)

\item Can a semi-disc be the numerical range of some matrix? 
That is: Let $A$ be an $n\times n$ complex matrix. Is it possible that the numerical range $W(A)=\{ x^*Ax \mid  x\in \Bbb C^n, x^*x=1\}\subset \Bbb C$ is a half disc?
($x^*$ stands for the transpose conjugate of the column vector $x$.) (See p.~197.)

\item
Are a square  matrix and its transpose  always similar? 
That is: Let $A$ be a square matrix over a field $\Bbb F$ and let $A^{\T}$ be the transpose of $A$. Does there exist an invertible
matrix $P$ (depending on $A$) over $\Bbb F$  such that $P^{-1}AP=A^{\T}$?
(See p.~203.)
\end{enumerate}

The book contains many nice, neat results in matrix theory. A lot of traditional topics and previously known results are treated in  new ways.
The contents of the book may be divided into four parts according to their topics:
\begin{itemize}
\item[I.] Chapters 1 and 2 provide the preliminaries and basics in matrix theory and linear algebra. They  prepare for the later chapters.
\item[II.] Chapters 3, 4, and 5 study matrix inequalities and lean to matrix analysis. Important matrix inequalities and their extensions of (mainly) majorization type are included.
\item[III.] Chapters 6, 7, and 8 are devoted to combinatorial aspects of matrices, including nonnegative matrices and graphs, matrix completion problems, and sign patterns.

\item[IV.] Chapters 9 and 10 present selected, miscellaneous topics, and show interaction of matrices with other fields of mathematics.
  Results are  interesting, diverse, even exciting.
\end{itemize}

Below are my recaps and highlights of each chapter in the book.
\medskip

{\bf Chapter 1} (Preliminaries) starts with basic concepts and notations for matrix theory, including special matrices,
 characteristic polynomials, norms, matrix decompositions, etc. Most, but not all,  of the materials are standard and can be found in  matrix or linear algebra texts. Some results (such as Gr\"{o}bner bases, see p.~25, p.~150) described in this chapter are used in later chapters, some are not; they are included here because they are interesting and their proofs are accessible to students. Take, for instance, Mirsky's result (on the existence of a matrix with prescribed diagonal entries and eigenvalues.
 Corollary 1.6, p.~9); its proof using Fillmore's theorem  (on similarity to a matrix with prescribed diagonal entries. Theorem 1.5, p.~8) appears to be new and is easily understood
 by students.
  Some topics and treatments may not be as familiar to the reader. For example, the idea of using topology to show the Cayley-Hamilton theorem is  clever and novel; while the notion of Gr\"{o}bner bases is not seen in standard matrix or linear algebra texts (but is found in \cite{Ada94}).

There is a short section about the books and journals in the field, directing the reader for further reading.

Problem 8  in the exercise section is the famous Gershgorin Disc Theorem (p.~32). We would like to
recommend \cite[p.~388]{HJ1.13} as a reference in which a rigourous  proof is provided. 
(Note that Ger\v{s}gorin or Gershgorin is spelled as Gergorin on p.~189 and p.~252 and as Gersgorin on p.~32 and p.~44.)
\medskip

{\bf Chapter 2} (Tensor Products and Compound Matrices) covers the tensor (a.k.a. Kronecker) product, Schur (a.k.a Hadamard) product,
 and compound matrices. These are important and essential tools in matrix theory and also in other fields such as quantum computation. In addition to basic properties,  the elegant Schur theorem (the Schur product of positive definite matrices is positive definite. Theorem 2.4, p.~38) is presented here, along with the fact  that a Schur product is a principal submatrix of a Kronecker product (Lemma 2.3, due to Marcus and Khan \cite{MarKhan59}). Through matrix ``vec",
it is shown that $A\otimes B$ and $B\otimes A$ are similar via a permutation matrix (which solely depends on $n$ when $A$ and $B$ are $n$-square matrices).
Included in this short chapter are also the Roth theorem (on the solvability of the matrix equation $AX-XB=C$. Theorem 2.8, p.~41), Lyapunov theorem (on the matrix equation $AX+XA^*=P$ with stable $A$ and positive definite $P$. Theorem 2.11, p.~43),
as well as Frobenius-K\"{o}nig theorem (on the size of a zero submatrix. Theorem 2.13, p.~46).
The Gershgorin  Disc Theorem (Corollary 9.11 of Chapter 9) is used in the proof of the Lyapunov theorem.

\medskip

{\bf Chapter 3} (Hermitian Matrices and Majorization) explores majorization theory, mainly for Hermitian matrices. Hermitian matrices are probably the most important family of matrices.
The chapter starts with the Courant-Fischer Min-Max expressions for the eigenvalues of Hermitian matrices, followed by the Cauchy Eigenvalue Interlacing Theorem (Theorem 3.2, p.~52)
and results of Weyl and Fan of similar type.

The  concept of majorization theory is introduced in Section 3.2. Its connection with doubly stochastic matrices is the theorem of Hardy-Littlewood-P\'{o}lya
(Theorem 3.9, p.~57). As the ($n\times n$) doubly stochastic matrices form a convex set, the author goes a bit further,  presenting a Krein-Milman result (on the extreme points of a compact convex set. Theorem 3.10) of
  convex analysis flavor.  The theorem of Birkhoff (every doubly stochastic matrix is a convex combination of permutation matrices. Theorem 3.11) is classic, and the proof is by induction on the number of zero entries. Combining  the Hardy-Littlewood-P\'{o}lya theorem and the Birkhoff theorem gives the Rado theorem (Theorem 3.12, p.~60)
   which states that if $x$ is majorized by $y$ then $x$ is a convex combination of the permutations of $y$. How many permutations of $y$? Theorem 3.13 says
   $n$ is enough (assuming that $x, y\in \Bbb R^n$) -- a nice result.

  The diagonal entries and eigenvalues of a square matrix are closely related. For Hermitian matrices, the diagonal entries are majorized by the eigenvalues
  (Schur theorem, Corollary 3.19, p.~62).
  Here is an inverse problem:
  Given two sets of real numbers, when can they be diagonal entries and eigenvalues of a matrix? Horn's theorem (Theorem 3.20, p.~63) gives a complete description in terms of majorization (with a proof by Chan-Li). To be exact, let $d=(d_1, \dots, d_n)$ and $\l=(\l_1, \dots, \l_n)$ be two sequences of real
  numbers. Then there exists an $n$-square real symmetric matrix with $d_1, \dots, d_n$ as diagonal entries and $\l_1, \dots, \l_n$ as eigenvalues if and only if $d$ is majorized by $\l$. This illustrates  the importance of majorization in matrix theory.

  The relation between convex function and majorization is revealed as Theorems 3.25 and 3.26 (pp.~65--66).
  Finally, the positive semidefinite matrices come into the play: the inequalities of L\"{o}wner-Heinz ($A^r\geq B^r$ for $A\geq B\geq 0$ and $0\leq r\leq 1.$ Theorem 3.29, p.~68) and Wang-Gong ($(\l(A^tB^t))^{1/t}$ is log-majorization increasing in $t$. Theorem 3.30, p.~69) on matrix (eigenvalue) monotonicity  are included in the last section of the chapter.

  Problem 11 in the Exercises (p.~74) suggests that the inequality of L\"{o}wner-Heinz (Theorem 3.29, p.~68) be proved by using an integral expression.
The reader should not underestimate the role of analysis in matrix theory. In fact, for many advanced matrix problems and results, analysis is an essential and efficient tool.

Problem 19 (p.~76) is the once long standing van der Waerden Conjecture (the permanent of an $n$-square doubly stochastic matrix is greater than or equal to $n!/n^n$. This was resolved by Egorychev and Falikman independently). The author guides the student to work on this problem
as a project and read more in that direction.

\medskip

{\bf Chapter 4} (Singular Values and Unitarily Invariant Norms) is a continuation of the study of inequalities (of eigenvalues for Hermitian matrices, Chapter 3) for singular values.
Many results of this chapter are singular value counterparts of eigenvalue inequalities.
  The results of Fan, Horn, and Weyl on singular values are presented in Section 4.1
 in terms of majorization.
The Gelfand theorem (Lemma 4.13, p.~84) on spectral radius and spectral norm of a matrix is stated and proved as a lemma to a result of Yamamoto (Theorem 4.14, p.~85) which shows a relation between the $i$th eigenvalue in absolute value and the
$i$th singular value (of the $k$th power as $k\rightarrow \infty$). The proof of Yamamoto's  theorem using compound matrices is clever and new.
The inequalities of Audenaert (on operator monotone functions and positive semidefinite matrices. Lemma 4.16, p.~86, Theorem 4.19, p.~87) are relatively new (not yet seen in other books).

Majorization inequalities, symmetric gauge functions, and unitarily invariant norms are closely related; roughly speaking, they are equivalent. These are discussed and shown by the von Neumann Theorem (a unitarily invariant norm is a symmetric gauge function acting on singular values. Theorem 4.23, p.~90) and the Ky Fan Dominance Principle ($\|A\|\leq \|B\|$ for all unitarily invariant norms if and only if it holds for Fan $k$-norms. Theorem 4.25, p.~92), followed by
the classic results of Mirsky (for any square matrices $A$ and $B$ of the same size and for any
unitarily invariant norm, $\|\diag (s(A)-s(B))\|\leq \|A-B\|$. Theorem 4.28, p.~94),
 Lidskii (for any $n\times n$ Hermitian matrices $G$ and $H$,
 $\l (G)-\l (H)\prec \l (G-H)$. Theorem 4.29, p.~95), and Toeplitz (for any matrix $A$, $w(A)\leq \|A\|_{\infty}\leq 2 w(A),$ where $w(A)$ is the
 numerical radius of $A$. Lemma~4.30, p.~95).

A matrix can be decomposed in various ways. The Cartesian decomposition is a decomposition
 that plays important roles in deriving matrix inequalities. This appears  in the results
 of Ando-Bhatia (on the eigenvalue and singular value majorization of   the real and imaginary parts in the  decomposition. Theorems 4.35, p.~97) and Bhatia-Kittaneh (on the Schatten $p$-norms of   the real and imaginary parts in the decomposition. Theorems 4.36, 4.37, pp.~98--99).

 There are many  interesting, elegant, and important inequalities selected in the Exercises (pp.~100--101).

\medskip

{\bf Chapter 5} (Perturbation of Matrices) deals with perturbation of matrices. A primary question is how ``far apart" or ``close" the spectra of two matrices can be in terms of matrix norm; one wants to estimate $d(\s(A), \s(B))=\min_{\mu \in \s(B)}\max_{\l \in \s(A)} |\l-\mu|$.
The chapter starts with a result of Bhatia-Elsner-Krause (Theorem 5.2, p.~104) which gives an upper bound of $d(\s(A), \s(B))$ in terms of the spectral norm. The proof of this theorem by using
a homotopic argument may be unfamiliar but it is interesting and valuable  to the reader. The central result of the chapter is the theorem of Hoffman-Wielandt (Theorem 5.3, p.~105). It states that
for any $n\times n$ normal matrices $A$ and $B$ with eigenvalues $\l_1, \dots, \l_n$ and $\mu_1, \dots, \mu_n$ (in any order), respectively,
$$\left ( \sum_{i=1}^n|\l_i-\mu_i|^2\right )^{1/2}\leq \|A-B\|_F.$$

Variations of the above result are discussed by Sun with $A$ being normal and $B$ arbitrary (Theorem 5.6, p.~108) and by Zhan with $T=A+iB$, where $A$ and $B$ are Hermitian (Theorem 5.9, p.~110), followed by a few interesting results of  Li (Theorem 5.10, p.~112) on the unitary matrices in the polar decomposition, Kittaneh (Theorem 5.12, p.~113), and Araki-Yamagami (Corollary 5.13, p.~113) on the positive semidefinite matrices in the polar  decomposition.
Norm estimations of band parts and backward perturbation analysis are also included in this chapter.
\medskip
\newpage

{\bf Chapter 6} (Nonnegative Matrices) studies matrices with real nonnegative entries; it includes the fundamental results of Perron, Frobenious, and Wielandt. Section 6.1 is devoted to
the Perron-Frobenius theory on nonnegative matrices. The heart of the theory is the Perron-Frobenius Theorem (Theorem 6.8, p.~123): If $A$ is an irreducible nonnegative matrix of order $n\geq 2$, then the spectral radius $\rho (A)>0$ is a simple eigenvalue of $A$ (called Perron root of $A$), and it has a positive corresponding eigenvector.
 The treatment of the material is standard and it is as elegant as it has always been. At the end of this section, the questions 1 and 2 we asked previously are answered:
 Every complex number is an eigenvalue of some nonnegative matrix (a proof due to Shan, p.~131); the square of a real matrix of order $n$ can have at most $n^2-1$ negative entries, and this upper bound can be attained (Theorem 6.24, due to Eschenbach-Li, p.~131).

 The relationship between  nonnegative matrices and digraphs is shown: A square (nonnegative) matrix is irreducible if and only if its digraph is strongly connected (Lemma 6.25, p.~133).
 Primitivity of nonnegative matrix is thoroughly discussed. A beautiful result of Wielandt states that if $A$ is primitive of order $n$, then
 $A^{(n-1)^2+1}>0$ (entrywise).  A collection of results about some special classes of nonnegative matrices such as totally nonnegative matrices, $M$-matrices, and $Z$-matrices are presented. It is notable that at the end of this chapter infinite divisibility of nonnegative symmetric matrices and complete monotonicity of functions and their roles in preserving positive semidefiniteness (Micchelli's result, Lemma 6.40, p.~146,  and Bapat's result, Theorem 6.41, p.~146) are briefly described. The Hadamard (entrywise) power (p.~145) of a nonnegative positive semidefinite matrix is extensively investigated in
 \cite{FitHor77}. Renewed interests in this direction are seen recently in
 \cite{StanfordAIM14}.
\medskip

{\bf Chapter 7} (Completion of Partial Matrices) is about completing a matrix in which some entries are known or prescribed and the matrix to be completed posses certain desired properties. For instance,
if the off-diagonal entries of a square matrix are known, and if we are allowed to choose freely  the diagonal entries, what can we say about the eigenvalues? Friedland's Theorem
(Theorem 7.1, p.~150) says that the eigenvalues can be anything. To be precise, let $\l_1, \dots, \l_n$ be elements of an algebraically closed field $\Bbb F$,  let $A$ be an $n\times n$ matrix over $\Bbb F$ whose off-diagonal entries are prescribed. Then one can always choose elements $a_1, \dots, a_n$ from the field as diagonal entries of the matrix $A$ so that $\l_1, \dots, \l_n$ are the eigenvalues of $A$. The proof of Friedland's key Lemma 7.3 (p.~150) by using Gr\"{o}bner bases is new to the reviewer.
This statement is shown through a few lemmas whose proofs use Buchberger's Theorems (Theorems 1.28 and 1.29, pp.~26--27) as well as Hilbert's Nullstellensatz (Theorem 10.6, p.~218).
Friedland's result is about diagonal completion, while Farahat-Ledermann's (Theorem 7.5, p.~153)  is about borderline completion with a prescribed characteristic polynomial.
The last section of the chapter is devoted to positive semidefinite completion,  which has much to do with graph theory. The Grone-Johnson-S\'{a}-Wolkowicz Theorem
(Theorem 7.21, p.~163) is a good example of results of this type.
\medskip

{\bf Chapter 8} (Sign Patterns) addresses  real matrices with (some) entries whose signs are known. What properties can one deduce for the (class of) matrix if the signs of some entries are known?  For instance, the diagonal entries of positive definite matrices have + signs
(so if there is a negative diagonal entry, then the matrix cannot be positive definite). A real matrix with sign pattern $\left [ {0 \atop -} {+\atop 0} \right ]$
has to be nonsingular and diagonalizable since its determinant is nonzero and the trace is zero (thus it has two different eigenvalues).
 A $6\times 6$ real matrix with prescribed sign pattern
but  no real eigenvalue(s) is shown (p.~165). This chapter overviews many aspects and results on the topic.
 A related, but not yet well understood problem is what sign patterns allow or require all distinct eigenvalues. (See Exercise (4), p.~179.)

\medskip

{\bf Chapter 9} (Miscellaneous Topics) presents several diverse matrix problems. The first section answers our Question 3  in the affirmative:
If two real matrices are (unitarily) similar over the complexes, then they are (real orthogonally) similar over the reals.
The commutator $[X, Y]=XY-YX$ of two matrices is studied and a few inequalities are shown.
The L\'{e}vy-Desplanques diagonal dominance theorem is discussed. The theorem of Camion-Hoffman (Theorem 9.14, p.~190) gives a   converse.
A simple version of the Gershgorin  Disc Theorem (Corollary 9.11, p.~189) is given as a consequence of the  L\'{e}vy-Desplanques Theorem.

Numerical range is a classic, but still active, research area in matrix theory. Section 9.5 studies  the shape of a numerical range.
When is a numerical range a line segment? (p.~193) When does an eigenvalue of a matrix lie on the boundary of its numerical range? (p.~195)
Can a semi-disc be the numerical range of some matrix? (Our Question 4. The answer is no; p.~197.)  These questions are thoroughly studied and answered.

Canonical forms of matrices are of central importance in matrix theory and linear algebra. The Smith canonical form (via invariant factors;
Theorem 9.26,
p.~200), rational canonical form (via companion matrix, Theorem 9.32, p.~204), and Jordan canonical form (via Jordan block; Theorem 9.34, p.~205) are studied.
As side-products, two results of Voss  are immediate:
{every square matrix over a field and its transpose are similar} (Corollary 9.30, p.~203. This answers our Question 5),
and {every square matrix over a field is a product of two symmetric matrices} (Theorem 9.35, p.~206). A similar question is:
Is every square matrix congruent to its transpose? The answer is yes. That is, if $A$ is a square matrix over a field, then there exists an invertible matrix $Q$ over the field such that
$Q^{\T}AQ=A^t$; see
\cite{DK02}. What about the $*$-congruence, i.e.,
$A^{\T}=R^*AR$ for some nonsingular $R$? The answer is also yes; see \cite{HS04}
for detail.

There are ``a lot" of zero entries in the aforementioned canonical forms.
At the end of the chapter, a question about the number of off-diagonal zero entries is asked and answered: under similarity, how many off-diagonal zero entries can a square matrix have at most? To be precise, let $A$ be a square matrix
over an algebraically closed field $\Bbb F$, and let $B$ be a matrix similar to $A$ over $\Bbb F$. Then the number of off-diagonal zero entries of $B$ is less than or equal to the number of off-diagonal zeros
in the Jordan canonical form of $A$. This is a nice result of Brualdi-Pei-Zhan (Theorem 9.38, p.~209).

The treatments of most materials in Chapter 9 (except for Section 9.7)  seem innovative.  
The author  compiles Sections 9.4 and 9.5 in a more readable way, though the ideas can be traced back to the original papers.

\medskip

{\bf Chapter 10} (Applications of Matrices), the last chapter,  showcases applications and connections of
matrix theory with combinatorics, number theory, algebra, and geometry.
This is not found in other standard textbooks.

 A combinatorial result of Erd\"{o}s-R\'{e}nyi-S\'{o}s (the Friendship Theorem. Theorem 10.1, p.~214)
is proved through the adjacency matrix of a graph;
the matrix proof of Theorem 10.2 (pp.~215--216) on counting elements in the intersections of $n$-sets is much shorter and more elementary than other proofs;
 the statement that the algebraic numbers form a field (Theorem 10.3, p.~216)
is proved through companion matrices and Kronecker products; the resultant of two polynomials is computed through a determinant, which is used in proving the Hilbert Nullstellensatz (Theorem 10.6, p.~218); in geometry, the volume of a polytope spanned by vectors (Theorem 10.9, p.~220)
 may be computed through determinants; a result of Cauchy (Theorem 10.13, p.~224)
on roots of a polynomial is proved by using  a nonnegative matrix whose digraph is strongly connected; and a result of Abel on a pair of complex polynomials (Theorem 10.15, p.~225)
is proved using   the Vandermonde matrix, while the original proof of Abel uses integrals.

\medskip

Finally, there is a section ({\bf Unsolved Problems}) as an appendix to the chapters. It contains a selection of unsolved problems; they are named as
conjectures, questions, and problems and numbered in a series. The very first one, Conjecture 1, is about the existence of the Hadamard matrices. This open question is very
important and well-known to mathematicians in combinatorics and geometry. It states that for every positive integer $n$
there exists a $4n$-square matrix of 1's and $-1$'s with rows (and columns) mutually orthogonal (a Hadamard matrix).
 A related question is: Can the permanent of a Hadamard matrix of order $n$ vanish for $n>2$? (Question 9, p.~231).
(Wanless has shown that the permanent of an $n\times n$ Hadamard matrix is nonzero if $2<n<32$; p.~231.)

Problem 3 asks for necessary and sufficient conditions that  a set of complex numbers  be the eigenvalues of a nonnegative matrix.

The ``permanental dominance conjecture" is presented as
Conjecture 6. This is also a long--open problem. A great amount of research work has been done;
 partial results have been obtained. The conjecture   is still open.

The Grone-Merris conjecture on Laplacian spectra (Conjecture 11, p.~231) is a problem involving matrices, majorization, combinatorics,
 and graph theory.
This intriguing problem has been resolved by H. Bai \cite{BaiTran}.

Regarding the distance of the spectra of two matrices with respect to norms, a problem of Bhatia (Problem 19, p.~234) is stated as:
Determine the smallest constant $c$ such that $d(\s(A), \s(B))\leq c\|A-B\|$ for normal matrices $A$ and $B$ of any order.

Exercises are assigned for each of Chapters 1--6 and Chapter 8; they are of various difficulties and depth; detailed literature,
 history and status, and recent developments are provided for most of the problems.

Matrix theory is a vast subject. It is impossible for a book of 250 pages to cover everything in the field.
This book is a collection of many important and interesting topics on matrices.
Students (and researchers as well) using the book will  gain and broaden 
their
  knowledge and skills in matrix theory and linear algebra.

Clearly, the author has followed his standards and criteria for this book: the materials included are  important,
 elegant,  ingenious, and  interesting. The contents of the book show great enthusiasm and excitement of the author for the subject.

In summary, the book provides an extensive survey and comprehensive treatment of matrix theory and its applications in a broad spectrum of fields.
I highly recommend it to the reader.

The reviewer is grateful to  Prof.  Roger Horn for his valuable comments.


\bigskip
\begin{flushright}
 Fuzhen Zhang\\
College of Arts and Sciences\\
 Nova Southeastern University\\
Fort Lauderdale, Florida 33314, USA\\
zhang@nova.edu\\
\end{flushright}
\end{document}